\documentclass[11pt]{amsart}
\usepackage{amssymb}
\usepackage{amsmath}

\newcommand{\mpn}{\medskip\par\noindent}
\newcommand{\pn}{\par\noindent}

\newcommand{\bpn}{\bigskip\par\noindent}
\theoremstyle{definition}

\theoremstyle{remark}

\numberwithin{equation}{section}

\begin{document}
\newcommand{\Mod}[1]{\,(\text{\mbox{\rm mod}}\;#1)}
\title[Umbral calculus and Euler polynomials]
{Umbral calculus and Euler polynomials}
\author{Dae San Kim, Taekyun KIM and Seog-Hoon Rim}
\begin{abstract}
In this paper, we study some properties of Euler polynomials arising
from umbral calculus. Finally, we give some interesting identities
of Euler polynomials using our results. Recently, Dere and Simsek
have studied umbral calculus related to special polynomials
$(see[6])$.
 \end{abstract}

 \maketitle

\section{\bf Introduction}
We recall that the Euler polynomials are defined by the generating function to be
\begin{align*}\tag{1}
\frac{2}{e^t+1}e^{xt}=e^{E(x)t}=\sum_{n=0}^{\infty}E_{n}(x)\frac{t^{n}}{n!}, \quad (see
[2,3,9-18]).
\end{align*}
with the usual convention about replacing $E^n(x)$ by $E_n(x)$.
\\In the special case, $x=0$, $E_{n}(0)=E_{n}$ are called the $n$-th
Euler numbers. From (1), we have
\begin{align*}\tag{2}
E_n(x)=\sum_{l=0}^{n} \binom{n}{l}x^{l}E_{n-l}=(E+x)^n, \quad
(see [4,5,7,8]).
\end{align*}
Thus, by (1) and (2), we get
\begin{align*}\tag{3}
E_0=1,\quad (E+1)^n+E_n=E_n(1)+E_n=2\delta_{0,n},
\end{align*}
where $\delta_{k,n} $ is the Kronecker's symbol $(see [2,3,8])$.
Note that $E_n(x)$ is a monic polynomials of degree $n$ and $E_n'(x)=\frac{dE_n(x)}{dx}=nE_{n-1}(x)$.
\\Let $\mathbb{C}$ be the complex number field and let $\mathcal{F}$ be the set of all formal power series in the variable $t$ over $\mathbb{C}$ with
\begin{align*}\tag{4}
\mathcal{F}=\{f(t)=\sum_{n=0}^{\infty}\frac{a_k}{k!}t^k ~|~a_k \in\mathbb{C}\}.
\end{align*}
Let $\mathbb{P}=\mathbb{C}[t]$ and let $\mathbb{P^*}$ be the vector
space of all linear functionals on $\mathbb{P}$. Now we use the notation 
$<L ~|~ p(x)>$ to denote the action of a linear functional $L$ on a polynomial
$p(x)$ $(see [6,13])$. We remind that the vector space operations on
$\mathbb{P^*}$ are defined by $<L+M ~|~ p(x)> = <L ~|~ p(x)> + <M
~|~ p(x)>$, $<cL ~|~ p(x)> = c<L ~|~ p(x)>$, $(see[6,13])$, where
$c$ is any constant in $\mathbb{C}$. The formal power series
\begin{align*}\tag{5}
f(t)=\sum_{k=0}^{\infty}\frac{a_k}{k!}t^k \in \mathcal{F},\quad (see [6,13]),
\end{align*}
defines a linear functional on $\mathbb{P}$ by setting
\begin{align*}\tag{6}
<f(t) ~ | ~x^n>=a_n ,\quad for\quad all\quad  n\geq 0.
\end{align*}
Thus, by (5) and (6), we have
\begin{align*}\tag{7}
<t^k ~ | ~x^n>=n!\delta_{n,k},  \quad (see [6,13]).
\end{align*}
Let $f_L(t)=\sum_{k=0}^{\infty}\frac{<L~|~ x^k>}{k!}t^k $. Then,
from (5) and (7), we note that \\$<f_L(t) ~|~ x^n> = <L~|~x^n>$ and
so as linear functionals $L=f_L(t)$ $(see[13])$. It is known in
[6,13] that the map $L \mapsto f_L(t)$ is a vector space isomorphism
from $\mathbb{P^*}$ onto  $\mathcal{F}$. Henceforth, $\mathcal{F}$
will denote both the algebra of formal power series in $t$ and the
vector space of all linear functionals on $\mathbb{P}$ and so an
element $f(t)$ of $\mathcal{F}$ will be thought of as both a formal
power series and a linear functional. We shall call $\mathcal{F}$
the umbral algebra $(see[6,13])$. The umbral calculus is the study
of umbral algebra and modern classical umbral calculus can be
described as a systematic study of the class of Sheffer sequences. By
(6), we see that
$$
<e^{yt}~|~x^n>=y^n \quad and \quad so \quad
 <e^{yt}~|~p(x)>=p(y) \quad (see [6,13]).
$$
Note that for all $f(t)$ in $\mathcal{F}$
\begin{align*}\tag{8}
f(t)=\sum_{k=0}^{\infty}\frac{<f(t)~|~x^k>}{k!}t^k ,
\end{align*}
and for all polynomials $p(x)$
\begin{align*}\tag{9}
p(x)=\sum_{k=0}^{\infty}\frac{<t^k~|~p(x)>}{k!}x^k ,\quad (see [6,13]).
\end{align*}
For $f_1(t),f_2(t), ... f_m(t) \in \mathcal{F}$, we have
\begin{align*}\tag{10}
 &< f_1(t)\times f_2(t)\times ...\times f_m(t)~|~x^n >
 \\& =\sum_{}^{ }\binom{n}{i_1,i_2,...i_m}<f_1(t)~|~x^{i_1}>...<f_m(t) ~|~x^{i_m}>, 
\end{align*}
where the sum is over all non-negative integers $i_1,i_2,...i_m$ such that $i_1+...+i_m=n$
$(see[6,13])$. The order $O(f(t))$ of the power series $f(t)\neq
0$ is the smallest integer $k$ for which $a_k$ does not vanish.
\\Now we set $O(f(t))=\infty$  if  $f(t)=0$. From the definition of order,we note that $O(f(t)g(t))=O(f(t)+O(g(t))$
and $O(f(t)+g(t))\geq min \{O(f(t), O(g(t))\}$. The series $f(t)$ has
a multiplicative inverse, denoted by $f(t)^{-1}$ or
$\frac{1}{f(t)}$, if and only if $O(f(t))=0$. Such a series is called an
invertible series. A series $f(t)$ for which $O(f(t))=1$ is
called a delta series $(see [6,13])$. let $f(t),g(t) \in
\mathcal{F}$. Then, we easily that $$<f(t)g(t) ~|~ p(x)> = <f(t)
~|~ g(t)p(x)> = <g(t) ~|~ f(t)p(x)>.$$ By (9), we set
\begin{align*}\tag{11}
p^{(k)}(x)= \frac{d^kp(x)}{dx^k}&=\sum_{l=k}^{\infty}\frac{1}{l!}<t^l~|~p(x)>l(l-1)...(l-k+1)x^{l-k}
\\&=\sum_{l=k}^{\infty}\frac{1}{(l-k)!}<t^k~|~p(x)>x^{l-k}.
\end{align*}
Thus, from (11), we have
\begin{align*}\tag{12}
p^{(k)}(0)= <t^k~|~p(x)> \quad and <1 ~|~ p^{(k)}(x)>=p^{(k)}(0).
\end{align*}
By (12), we get
\begin{align*}\tag{13}
t^kp(x)=p^{(k)}(x)= \frac{d^kp(x)}{dx^k}\quad (see [6,13]).
\end{align*}
Thus, by (13) and Taylor expansion, we get
\begin{align*}\tag{14}
e^{yt}p(x)=p(x+y)\quad (see [6,13]).
\end{align*}
Let $S_n(x)$ be polynomials with deg $S_n(x)=n$ and let $f(t)$ be a
delta series and $g(t)$ be an invertible series. Then there exists a
unique sequence $S_n(x)$ of polynomials with $< g(t) f(t)^k
~|~S_n(x)>=n!\delta _{n,k}, (n,k \geq 0)$. The sequence $S_n(x)$ is
called the Sheffer sequence for $(g(t),f(t))$, which is denoted by
$S_n(x)\sim (g(t),f(t))$ $(see[13])$. If $S_n(x) \sim (1,f(t))$,
then $S_n(x)$ is called the associated sequence for $f(t)$, or
$S_n(x)$ is associated to $f(t)$. If $S_n(x) \sim (g(t), t)$, then
$S_n(x)$ is called the Appell sequence for $g(t)$ or $S_n(x)$ is
Applell for $g(t)$ $(see [13])$. For $p(x) \in \mathbb{P}$, it is
known in [13] that
\begin{align*}\tag{15}
<\frac{e^{yt}-1}{t} ~|~p(x)>= \int_0^{y}p(u)du,
\end{align*}
\begin{align*}\tag{16}
<f(t) ~|~xp(x)>= <\partial_tf(t) ~|~ p(x)>=<f'(t) ~|~p(x)>,
\end{align*}
and
\begin{align*}\tag{17}
<f(t) ~|~p(\alpha x)>= <f(\alpha t) ~|~ p(x)>,<e^{yt}-1~|~p(x)>=p(y)-p(0),
\end{align*}
where $\alpha$ is a complex constant $(see [6,13])$.
\\Suppose that $S_n(x)$ is the Sheffer sequence for $(g(t),f(t))$. Then we have the following equations:
\begin{align*}\tag{18}
h(t)=\sum_{k=0}^{\infty}\frac{<h(t)~|~S_k(x)>}{k!}g(t)f(t)^k ,\quad h(t)\in \mathcal{F},
\end{align*}
\begin{align*}\tag{19}
p(x)=\sum_{k=0}^{\infty}\frac{<g(t)f(t)^k~|~p(x)>}{k!}S_k(x) ,\quad p(x)\in \mathbb{P},
\end{align*}
\begin{align*}\tag{20}
f(t)S_n(x)=nS_{n-1}(x),
\end{align*}
and
\begin{align*}\tag{21}
\frac{1}{g(\overline{f}(t))}e^{y \overline{f}(t)}=\sum_{k=0}^{\infty}\frac{S_k(y)}{k!}t^k ,\quad for \quad all \quad  y\in \mathbb{C},
\end{align*}
where $ \overline{f}(t)$ is the compositional inverse of $f(t)$, $(see[13])$.
\\Recently, Dere and Simsek have studied umbral calculus related to special polynomials $(see[6])$.
In this paper, we study properties of Euler polynomials arising from umbral calculus. From our investigation, we derive some interesting identities of Euler polynomials.

\section{\bf  Umbral calculus and Euler polynomials}
Let $S_n(x)$ be an Appell sequence for $g(t)$. From (21), we note that
\begin{align*}\tag{22}
\frac{1}{g(t)}x^n=S_n(x) \Leftrightarrow x^n=g(t)S_n(x),\quad (n\geq 0).
\end{align*}
Let us take $g(t)=\frac{e^t+1}{2} \in \mathcal{F}$. Then, we note that $g(t)$ is an invertible series.
By (1), we get
 \begin{align*}\tag{23}
\sum_{k=0}^{\infty}E_k(x)\frac{t^k}{k!}=\frac{1}{g(t)}e^{xt},
\end{align*}
and, from (23), we have
\begin{align*}\tag{24}
\frac{1}{g(t)}x^n=E_n(x),(n\geq 0),\quad
tE_n(x)=n\frac{1}{g(t)}x^{n-1}=nE_{n-1}(x).
\end{align*}
Thus, by (21) and (24), we see that
$E_n(x)$ is the Appell sequence for $(\frac{e^t+1}{2},t)$.
Indeed,
\begin{align*}\tag{25}
<\frac{e^t+1}{2}t^k ~|~ E_n(x) >&=\frac{k!\binom{n}{k}}{2}<e^t+1 ~|~E_{n-k}(x)>
\\&=\frac{k!\binom{n}{k}}{2}(E_{n-k}(1)+E_{n-k}), \quad (n,k \geq 0).
\end{align*}
By (3) and (25), we get
\begin{align*}\tag{26}
<\frac{1+e^t}{2}t^k ~|~ E_n(x) >=n! \delta _{n,k}, \quad (n,k \geq
0).
\end{align*}
Therefore, by  (24) and (26), we obtain the following lemma.
\par\bigskip
{\bf Lemma  1 }. For $n \geq 0$ ,  $E_n(x)$ is the Appell sequence  for $(\frac{1+e^t}{2},t)$.
\par\bigskip
{\bf Remark }. By (1), we get
\begin{align*}\tag{27}
<\frac{2}{e^t+1} ~|~ x^n > =\sum_{l=0}^{n}\frac{E_l}{l!} <t^l~| ~x^n>=\sum_{l=0}^{n}\frac{E_n}{l!}n!\delta_{n,l}.
\end{align*}
Thus, by (27), we have
\begin{align*}
<\frac{2}{e^t+1} ~|~ x^n > =E_n, \quad (n\geq 0) .
\end{align*}
From Lemma l, we note that
\begin{align*}\tag{28}
\sum_{k=0}^{\infty}\frac{E_k(x)}{k!} t^k= \frac{1}{g(t)}e^{xt},
\end{align*}
where $g(t)=\frac{e^t+1}{2} \in \mathcal{F}$.
\\Let us take the first derivative with respect to $t$ on both sides in (28). Then we have
\begin{align*}\tag{29}
\sum_{k=1}^{\infty}\frac{E_k(x)}{k!} kt^{k-1}&= \frac{xg(t)e^{xt}-g'(t)e^{xt}}{g(t)^2}
\\&= \sum_{k=0}^{\infty}\{x\frac{1}{g(t)}x^k-\frac{g'(t)}{g(t)}\frac{1}{g(t)}x^k\}\frac{t^k}{k!}.
\end{align*}
Thus, by (28) and (29), we get
\begin{align*}
E_{k+1}(x)=\left(x-\frac{g'(t)}{g(t)}\right)E_k(x), \quad (k\geq 0).
\end{align*}
From (2) and (24), we get
\begin{align*}\tag{30}
\int_x^{x+y}E_n(u)du &= \frac{1}{n+1}\{E_{n+1}(x+y)-E_{n+1}(x)\}
\\&=\frac{1}{n+1}\sum_{k=1}^{\infty}\binom{n+1}{k}E_{n+1-k}(x)y^k
\\&=\sum_{k=1}^{\infty}(n)_{k-1}E_{n+1-k}(x)\frac{y^k}{k!}
\\&=\frac{1}{t} \left(\sum_{k=0}^{\infty} \frac{y^k}{k!}t^k-1 \right) E_n(x)=\frac{e^{yt}-1}{t}E_n(x).
\end{align*}
Therefore, by (30), we obtain the following lemma.
\par\bigskip
{\bf Lemma  2 }. For $n \geq 0$ , we have
\begin{equation*}
\begin{split}
\int_x^{x+y} E_n(u)du=\frac{e^{yt}-1}{t}E_n(x).
\end{split}
\end{equation*}
By (24), we get
\begin{align*}\tag{31}
E_n(x)=t\frac{1}{n+1}E_{n+1}(x).
\end{align*}
Thus, by (31), we get
\begin{align*}\tag{32}
&<\frac{e^{yt}-1}{t} ~|~ E_n(x) >=<\frac{e^{yt}-1}{t}|~ t\frac{1}{n+1}E_{n+1}(x)>
\\&=<e^{yt}-1~|~\frac{E_{n+1}(x)}{n+1}>=\frac{1}{n+1} \{ E_{n+1}(y)-E_{n+1}(0) \}
\\&=\int_0^yE_n(u)du.
\end{align*}
Therefore, by (30) and (32), we obtain the following proposition.
\par\bigskip
{\bf  Proposition 3 }. For $n \geq 0$ , we have
\begin{equation*}
\begin{split}
\int_{0}^{y}E_n(u)du=< \frac{e^{yt}-1}{t}~|~E_n(x)>.
\end{split}
\end{equation*}
Let
\begin{align*}\tag{33}
\mathbb{P}_n=\{p(x) \in \mathbb{Q}[x] ~|~ \deg p(x) \leq n\}.
\end{align*}
Then $\mathbb{P}_n$ is a $(n+1)$-dimensional vector space over $\mathbb{Q}$ and
 $\{E_0(x),E_1(x),..E_n(x)\}$ is a basis for $\mathbb{P}_n$.
 For $ p(x) \in \mathbb{P}_n$, we write it as
\begin{align*}\tag{34}
p(x) = \sum_{k=0}^{n}b_kE_k(x).
\end{align*}
By Lemma 1, we get
\begin{align*}\tag{35}
<\frac{e^t+1}{2}t^k ~|~ E_n(x) >=n! \delta _{n,k}, \quad (n,k \geq 0).
\end{align*}
From (34) and (35), we have
\begin{align*}\tag{36}
<\frac{e^t+1}{2}t^k ~|~ p(x) > &=\sum_{l=0}^{n}b_l <\frac{e^t+1}{2}t^k ~|~ E_l(x) > \\&=\sum_{l=0}^{n}b_l l!\delta_{l,k}=k! b_k.
\end{align*}
Thus, by (36), we get
\begin{align*}\tag{37}
b_k&=\frac{1}{k!}<\frac{e^t+1}{2}t^k ~|~ p(x) > =\frac{1}{2k!} <(e^t+1)t^k ~|~ p(x) >
\\&=\frac{1}{2k!}<e^t+1 ~|~ p^{(k)}(x)>=\frac{1}{2k!}\{p^{(k)}(1)+p^{(k)}(0)\}.
\end{align*}
Therefore, by (37), we obtain the following theorem.
\par\bigskip
{\bf Theorem 4 }. For $p(x) \in \mathbb{P}_n$ , let
\\
\begin{equation*}
\begin{split}
p(x) = \sum_{k=0}^{n}b_kE_k(x).
\end{split}
\end{equation*}
Then, we have
\begin{equation*}
\begin{split}
b_k=\frac{1}{2k!} <(e^t+1)t^k ~|~ p(x) > .
\end{split}
\end{equation*}
In other words,
\begin{equation*}
\begin{split}
b_k=\frac{1}{2k!}\{p^{(k)}(1)+p^{(k)}(0)\}.
\end{split}
\end{equation*}
\par\bigskip
Let us take $p(x)=B_n(x) \in \mathbb{P}_n$. Then $B_n(x)$ can be written as a linear combination of
$E_0(x), E_1(x), .. E_n(x)$ as follows:
\begin{align*}\tag{38}
B_n(x)= \sum_{k=0}^{n}b_kE_k(x),\quad (n\geq 0).
\end{align*}
By Theorem 4, we get
\begin{align*}
b_k&=\frac{1}{2k!}<(e^t+1)t^k ~|~ B_n(x) >
 \\&=\frac{1}{2k!}<e^t+1~|~ n(n-1)\cdots(n-k+1)B_{n-k}(x) >
\\&=\frac{1}{2k!}\binom{n}{k}k! <e^t+1~|~ B_{n-k}(x)>.
\end{align*}
From (17) and (38), we have
\begin{align*}\tag {39}
b_k=\frac{1}{2k!}\binom{n}{k}k! \{B_{n-k}(1)+B_{n-k} \}.
\end{align*}
 As is well known, the recurrence of
Bernoulli numbers is given by
\begin{align*}\tag {40}
B_0=1, B_n(1) -B_n=\delta_{1,n}, \quad (see[1,8]).
\end{align*}
By (38),(39) and (40), we get
\begin{align*}\tag {41}
B_n(x)&= b_nE_n(x)+b_{n-1}E_{n-1}(x)+\sum_{k=0}^{n-2}b_kE_k(x)
\\&=E_n(x)+ \sum_{k=0}^{n-2}\binom{n}{k} B_{n-k}E_k(x).
\end{align*}
Therefore, we obtain the following corollary.
\par\bigskip
{\bf Corollary 5 }. For $n \geq 0$, we have
\\
\begin{equation*}
\begin{split}
B_n(x)=E_n(x)+ \sum_{k=0}^{n-2}\binom{n}{k} B_{n-k}E_k(x).
\end{split}
\end{equation*}
\par\bigskip
For $r \in \mathbb{Z}_{+}=\mathbb{N} \cup \{0\}$, the Euler polynomials $E_n^{(r)}(x)$ of order $r$ are defined by the generating function to be
\begin{align*}\tag {42}
\underbrace{\left( \frac{2}{e^t+1}\right) \times \left( \frac{2}{e^t+1}\right)\times ... \times\left( \frac{2}{e^t+1}\right)}_{r-times}e^{xt}= \sum_{n=0}^{\infty}E_n^{(r)}(x)\frac{t^n}{n!}.
\end{align*}
In the special case, $x=0$, $E_n^{(r)}(0)=E_n^{(r)}$ are called the
$n$-th Euler numbers of order $r$ $(see [9,10])$. Let us take
\begin{align*}\tag {43}
g^r(t)=\frac{1}{2^r}\underbrace{(e^t+1) \times  ... \times (e^t+1)}_{r-times}= \left( \frac{e^t+1}{2}\right)^r.
\end{align*}
Then, we note that $g^r(t)$ is an invertible series.
\\From (42) and (43), we note that
\begin{align*}\tag {44}
\sum_{n=0}^{\infty}E_n^{(r)}(x)\frac{t^n}{n!}
=\frac{1}{g^r(t)}e^{xt}= \sum_{n=0}^{\infty}\frac{1}{g^r(t)}x^n\frac{t^n}{n!}.
\end{align*}
By (44), we get
\begin{align*}\tag {45}
E_n^{(r)}(x)=\frac{1}{g^r(t)}x^n, \quad tE_n^{(r)}(x)=\frac{n}{g^r(t)}x^{n-1}=nE_{n-1}^{(r)}(x).
\end{align*}
From (21) and (45), we note that $E_n^{(r)}(x)$ is the Appell sequence for $\left(\frac{e^t+1}{2}\right)^r$.
By Appell identity, we also get
\begin{align*}\tag {46}
E_n^{(r)}(x+y)=\sum_{k=0}^{n}\binom{n}{k}E_{n-k}^{(r)}(x)y^k.
\end{align*}
It is easy to show that
\begin{align*}\tag {47}
< \underbrace{\frac{2}{e^t+1}\times ...\times \frac{2}{e^t+1}}_{r-times}~|~x^n>=<\left(\frac{2}{e^t+1}\right)^r~|~x^n>=E_n^{(r)}.
\end{align*}
By (47), we get
\begin{align*}\tag {48}
<\left(\frac{2}{e^t+1}\right)^re^{yt}~|~x^n>=\sum_{l=0}^{n}\binom{n}{l}y^lE_{n-l}^{(r)}=E_n^{(r)}(y), \quad (n\geq 0),
\end{align*}
and
\begin{align*}\tag {49}
&<\left(\frac{2}{e^t+1}\right)^re^{yt}~|~x^n> =<\left(\frac{2}{e^t+1}\right)^re^{(\frac{y}{r})rt}~|~x^n>
\\&=\sum_{i_1+\cdots+i_r=n}^{}\binom{n}{i_1,...,i_r}<\left(\frac{2}{e^t+1}\right)e^{(\frac{y}{r})t}~|~x^{i_1}>...<\left(\frac{2}{e^t+1}\right)e^{(\frac{y}{r})t}~|~x^{i_r}>
\\&=\sum_{i_1+\cdots+i_r=n}^{}\binom{n}{i_1,...,i_r}E_{i_1}(\frac{y}{r})...E_{i_r}(\frac{y}{r}).
\end{align*}
From (48) and (49), we have
\begin{align*}\tag {50}
E_n^{(r)}(x)=\sum_{i_1+\cdots+i_r=n}^{}\binom{n}{i_1,...,i_r}E_{i_1}(\frac{x}{r})...E_{i_r}(\frac{x}{r}).
\end{align*}
Therefore, we obtain the following theorem.
\par\bigskip
{\bf Theorem 6}. For $n \geq 0$,  $E_n^{(r)}(x)$ is the Appell sequence for $\left(\frac{e^t+1}{2}\right)^r$.
\\In addition,
\\
\begin{equation*}
\begin{split}
E_n^{(r)}(x)=\sum_{i_1+\cdots+i_r=n}^{}\binom{n}{i_1,...,i_r}E_{i_1}(\frac{x}{r})...E_{i_r}(\frac{x}{r}).
\end{split}
\end{equation*}
\par\bigskip
From (46) and (50), we note that $E_n^{(r)}(x)$ is a monic
polynomials of degree $n$.
\\Let $p(x)=E_n^{(r)}(x) \in \mathbb{P}_n$. Then $E_n^{(r)}(x)$ can be written as a linear combination of $E_0(x), ... ,E_n(x)$ as follows:
\begin{align*}\tag {51}
p(x)=E_n^{(r)}(x)=\sum_{k=0}^{n}b_kE_k(x), \quad (n\geq 0).
\end{align*}
By Theorem 4, $b_k$ is given by
\begin{align*}\tag {52}
b_k&=\frac{1}{2k!}<(e^t+1)t^k ~|~ p(x) > =\frac{1}{2k!} <(e^t+1)t^k~|~ E_n^{(r)}(x) >
\\&=\frac{\binom{n}{k}}{2}<e^t+1~|~ E_{n-k}^{(r)}(x) > =\frac{\binom{n}{k}}{2}\{E_{n-k}^{(r)}(1)+E_{n-k}^{(r)} \}.
\end{align*}
From (52), we have
\begin{align*}\tag {53}
\sum_{n=0}^{\infty} \{E_n^{(r)}(x+1)+E_n^{(r)}(x)\}\frac{t^n}{n!}=\left(\frac{2}{e^t+1}\right)^r(e^t+1)e^{xt}
\\=2\left( \frac{2}{e^t+1} \right)^{r-1}e^{xt}=2\sum_{n=0}^{\infty}E_n^{(r-1)}(x)\frac{t^n}{n!}.
\end{align*}
By comparing coefficients of both sides of (52), we get
\begin{align*}\tag {54}
E_n^{(r)}(x+1)+E_n^{(r)}(x)=2E_n^{(r-1)}(x), \quad (n\geq 0).
\end{align*}
Thus, from (51)and (54), we have
\begin{align*}\tag {55}
b_k=\binom{n}{k}E_{n-k}^{(r-1)}.
\end{align*}
Thus, by (51) and (55), we obtain the following corollary.
\par\bigskip
{\bf Corollary 7}. For $n \geq 0$, we have
\\
\begin{equation*}
\begin{split}
E_n^{(r)}(x)=\sum_{k=0}^{n}\binom{n}{k}E_{n-k}^{(r-1)}E_k(x).
\end{split}
\end{equation*}
\par\bigskip
It is not difficult to show that $\{ E_0^{(r)}(x), E_1^{(r)}(x),\cdots,E_n^{(r)}(x)\}$ is a basis for $\mathbb{P}_n$.
Let $p(x) \in \mathbb{P}_n$. Then $p(x)$ can be written as a linear combination of $ E_0^{(r)}(x), E_1^{(r)}(x),\cdots, E_n^{(r)}(x)$ as follows:
\begin{align*}\tag {56}
p(x)=\sum_{k=0}^{n}b_k^rE_k^{(r)}(x).
\end{align*}
By Theorem 6, we see that
\begin{align*}\tag {57}
< \underbrace{\frac{e^t+1}{2}\times ...\times \frac{e^t+1}{2}}_{r-times}t^k~|~E_n^{(r)}(x)>=n!\delta_{n,k},\quad (n,k\geq 0).
\end{align*}
From (56) and (57), we have
\begin{align*}\tag {58}
<\left(\frac{e^t+1}{2}\right)^r t^k ~|~p(x)> &=\sum_{l=0}^{n}b_l^r<\left(\frac{e^t+1}{2}\right)^r t^k ~|~E_l^{(r)}(x)>
\\&=\sum_{l=0}^{n}b_l^rl! \delta_{l,k}=k!b_k^r.
\end{align*}
By (58), we get
\begin{align*}\tag{59}
b_k^r&=\frac{1}{k!}<\underbrace{\frac{e^t+1}{2}\times... \times\frac{e^t+1}{2}}_{r-times} t^k ~|~ p(x) >
 \\&=\frac{1}{2^rk!} <(e^t+1)^rt^k ~|~ p(x) >
\\&=\frac{1}{2^rk!}\sum_{l=0}^{r}\binom{r}{l} <e^{lt} ~|~ t^kp(x)>
\\&=\frac{1}{2^rk!}\sum_{l=0}^{r} \binom{r}{l}<e^{lt} ~|~ p^{(k)}(x)>
\\&=\frac{1}{2^rk!}\sum_{l=0}^{r} \binom{r}{l}p^{(k)}(l),
\end{align*}
where $p^{(k)}(l)=\frac{d^kp(x)}{dx^k}|_{x=l}$.
Therefore, by (56) and (59), we obtain the following theorem.
\par\bigskip
{\bf Theorem 8}. For $p(x) \in \mathbb{P}_n$, let
\\
\begin{equation*}
\begin{split}
p(x) = \sum_{k=0}^{n}b_k^rE_k^{(r)}(x).
\end{split}
\end{equation*}
Then, we have
\begin{equation*}
\begin{split}
b_k^r=\frac{1}{2^rk!} <(e^t+1)^rt^k ~|~ p(x) > .
\end{split}
\end{equation*}
In other words,
\\$b_k^r=\frac{1}{2^rk!}\sum_{l=0}^{r} \binom{r}{l}p^{(k)}(l)$, where $p^{(k)}(l)=\frac{d^kp(x)}{dx^k}|_{x=l}$.
\par\bigskip
Let us take $p(x)=E_n(x) \in \mathbb{P}_n, \quad (n\geq 0).$
\\From Theorem 8, we have
\begin{align*}\tag {60}
E_n(x)=\sum_{k=0}^{n}b_k^rE_k^{(r)}(x),
\end{align*}
where
\begin{align*}\tag {61}
b_k^r&=\frac{1}{2^rk!}<(e^t+1)^rt^k ~|~ E_n(x)>
\\&=\frac{\binom{n}{k}}{2^r}<(e^t+1)^r ~|~ E_{n-k}(x)>
\\&= \frac{\binom{n}{k}}{2^r}\sum_{l=0}^{r} \binom{r}{l}E_{n-k}(l).
\end{align*}
By (60) and (61), we get
\begin{align*}
E_n(x)&=\frac{1}{2^r}\sum_{k=0}^n\binom{n}{k}\left( \sum_{l=0}^{r} \binom{r}{l}E_{n-k}(l)\right)E_k^{(r)}(x).
\end{align*}
Let $\alpha(\neq0) \in \mathbb{C}$. Then we have
\begin{align*}
E_n(\alpha x)=\alpha^n \frac{g(t)}{g(\frac{t}{\alpha})}E_n(x),
\end{align*}
where $g(t)=\frac{1}{2}(e^t+1).$
Let us consider the Bernoulli polynomials of order $s$ with
\begin{align*}
p(x)=B_n^{(s)}(x) \in \mathbb{P}_n, \quad (n\geq 0).
\end{align*}
By Theorem 8, we get
\begin{align*}\tag {62}
B_n^{(s)}(x)&=\sum_{k=0}^{n} b_k^rE_k^{(r)}(x),
\end{align*}
where
\begin{align*}\tag{63}
b_k^r&=\frac{1}{2^rk!} <(e^t+1)^rt^k ~|~ B_n^{(s)}(x) >
\\&=\frac{\binom{n}{k}}{2^r}<(e^t+1)^r ~|~ B_{n-k}^{(s)}(x)>
\\&=\frac{\binom{n}{k}}{2^r}\sum_{l=0}^{r}\binom{r}{l} <e^{lt} ~|~ B_{n-k}^{(s)}(x)>
\\&=\frac{\binom{n}{k}}{2^r}\sum_{l=0}^{r}\binom{r}{l}  B_{n-k}^{(s)}(l).
\end{align*}
Therefore, by (62) and (63), we get
\begin{align*}
B_n^{(s)}(x)=\frac{1}{2^r}\sum_{k=0}^{n}\binom{n}{k} \{\sum_{l=0}^{r}\binom{r}{l} B_{n-k}^{(s)}(l)\} E_k^{(r)}(x).
\end{align*}

\par\noindent
{\mpn { \bpn {\small Dae San {\sc Kim} \mpn Department of Mathematics,
\pn Sogang University, Seoul 121-742, S. Korea \pn {\it E-mail:}\
{\sf dskim@sogang.ac.kr} }
 \mpn { \bpn {\small Taekyun {\sc KIM} \mpn
 Department of Mathematics,\pn
Kwangwoon University, Seoul 139-701, S.Korea
  \pn {\it E-mail:}\ {\sf tkkim@kw.ac.kr} }
\mpn { \bpn {\small Seog-Hoon {\sc Rim} \mpn Department of
Mathematics Education, \pn Kyungpook National University, Taegu
702-701, S. Korea \pn {\it E-mail:}\ {\sf shrim@knu.ac.kr} }

\

\end{document}